\documentclass[11pt]{article}
\usepackage{amssymb}
\newcommand{\complex}{{\mathbb C}}

\newcommand{\integers}{{\mathbb Z}}

\newcommand{\delbar}{\overline{\partial}}
\newcommand{\indx}{{\rm index}}
\title{{\bf Some questions about the index of quantized contact
transformations\\}
}
\author{Alan Weinstein\thanks{This is an expanded version of a
lecture given at the Symposium on Geometric Methods in
Asymptotic Analysis, RIMS, Kyoto, May 20, 1997.
Research partially supported by NSF
Grant DMS-96-25122 and a JSPS Invitation Fellowship.  I would like to
thank RIMS (Kyoto University) and Keio University for their
hospitality.}
\\Department of Mathematics\\
University of California\\
Berkeley, CA 94720 USA\\
{\small(alanw@math.berkeley.edu)}}
\begin{document}
\maketitle
\section{Introduction}
\label{sec-intro}
If $M_1$ and $M_2$ are compact differentiable manifolds, a contact
diffeomorphism $\phi$ between their cosphere bundles gives rise to a
class $C(\phi)$ of Fredholm operators, called {\em Fourier integral
operators} or 
{\em quantized contact transformations} between the Hilbert spaces of
$L^2$ functions (or, more invariantly, half densities) on $M_1$ and
$M_2$.  The question of whether there is a unitary operator in this
class was raised in \cite{we:Fourier}, where such operators were used
to approximately intertwine the laplacians on riemannian manifolds
with symplectically equivalent geodesic flows.  It was shown there
that the existence of the unitary operator was equivalent to the
vanishing of the index of operators in $C(\phi)$, and the problem of
finding a topological formula for the index of the operators in
$C(\phi)$ was posed.  A conjecture for such a formula was made by
M.~Atiyah in a conversation with the author at some time in the
mid-1970's.  Little progress has been made since then, partly because
it is hard to produce examples where the index even has a chance
of being non-zero.

Recent developments in analysis and symplectic geometry have suggested
generalizations of this index problem to settings where non-zero
indices are known to exist, and technical advances in analysis seem to
have brought a solution within reach.  This talk will give an overview
of the problem and describe prospects for its solution in the context
of Epstein's relative index for CR structures \cite{ep:relative}.
Work of Guillemin \cite{gu:toeplitz} using analysis on Grauert tubes
implies that our original index problem can be set in this context.

Much of this paper is speculative in nature.  It is in part a report
on ongoing discussions (in person and by electronic mail) with David
Borthwick, Ana Cannas da Silva, Charles Epstein, Victor Guillemin,
and
Steven Zelditch.  I would like to thank all of them for their
contributions to this project.  In addition, I have received helpful
advice from Michael Christ, Peter Gilkey, Ian Grojnowski, Janos
Kollar, Richard Melrose, Gregory Sankaran, Bernard Shiffman, and
Sidney Webster.

\section{Polarizations of contact manifolds}

In this section, we will see how the index problem for Fourier
integral operators can be considered as a version of the question
``how does the quantum Hilbert space depend on the polarization?''
which is central to the theory of geometric quantization.  First of
all, we will recall how the notions of geometric quantization are
transplanted to contact manifolds from their usual symplectic setting.
This discussion is very much inspired by the work of Boutet de Monvel
and Guillemin \cite{bo-gu:spectral}. 

Let $Y$ be a contact manifold, $C\subset TY$ the contact distribution.
The bracket of sections of $C$ determines a natural nondegenerate
2-form $\Omega$ on $C$ with values in the normal line bundle $TY/C$.
A {\bf polarization} of $Y$ is defined to be a complex subbundle
${\cal J}$ of the complexification $C_{\complex}$ such that:
\begin{itemize}
\item (the natural complex extension of) $\Omega$ is zero on ${\cal J}$;
\item $\dim {\cal J} = \frac {1}{2} \dim C_{\complex}$;
\item $[\Gamma({\cal J} ),\Gamma({\cal J} )] \subseteq \Gamma({\cal J}
)$. 
\end{itemize}

One should add a further condition relating ${\cal J} $ and $ \overline
{\cal J} $, analogous to that in the symplectic case, but it will be
automatically satisfied in the two extreme cases which will interest
us in this paper.  

The ``quantum Hilbert space'' associated to the polarization ${\cal
J}$ is obtained by taking the space of smooth functions on $Y$ which
are annihilated by all sections of ${\cal J}$, and then taking its
closure $H_{\cal J} $ in  $L^2(Y)$ (defined with the aid of a chosen
volume element on $Y$).  A fundamental problem in geometric
quantization theory is to relate the Hilbert spaces arising from
different polarizations of the same contact manifold.  In our setting,
these spaces are infinite-dimensional, but we can define the
``difference between the dimensions'' of two such spaces as the index
of the orthogonal projection operator (in $L^2(Y)$) from one space to
the other.  We will call this index the {\bf relative index} of the
two polarizations.  We will see that, in many cases, the projection
operator is Fredholm, so that the relative
index is finite , and we will propose a topological formula for
computing it.  

Our basic idea is to associate to each polarization ${\cal J}_i$ of a
compact contact manifold $Y$ some ``filling'' of $Y$, i.e. some
compact manifold $X_i$ having $Y$ as its boundary.  The relative index
of two polarizations, defined provisionally as the index of the
orthogonal projection from one quantum Hilbert space to the other,
should then be the index of a Dirac operator on the manifold obtained
by gluing the two fillings along $Y$.  This is our {\bf gluing conjecture}.

\section{Complex polarizations}
\label{sec-complex}

A polarization ${\cal J}$ is called a {\bf complex polarization} if
${\cal J} $ and $ \overline {\cal J} $ are complementary subbundles.
Such polarizations are also known as (nondegenerate) {\bf CR} (or
Cauchy-Riemann) {\bf structures}.\footnote{For the most general CR
structures, $C\subset TY$ may be any distribution of codimension 1, not
necessarily contact.}  These complex polarizations almost complex
structures $J$ on the vector bundle $C$ by the rule ${\cal
J}=\{x-iJx|x \in C\}$.  The condition $[\Gamma({\cal J} ),\Gamma({\cal
J} )] \subseteq \Gamma({\cal J} )$ is the usual integrability
condition for CR structures.

For a complex polarization of CR type, the smooth functions
annihilated by the sections of ${\cal J}$ are generally known as CR
functions.  Their closure $H_J$ in $L^2(Y)$ is essentially independent of
the choice of volume element on $Y$ and is called the {\bf Hardy
space} of the CR structure.  The orthogonal projection onto this
quantum Hilbert space does depend on the volume element and is known as the
{\bf Szego projector}.

An important supplementary condition on complex polarizations is
strict {\bf pseudoconvexity}, which is definiteness of the
$TY/C$-valued {\bf Levi form} on $C$ defined by $(x,y)\mapsto
\Omega(Jx,y)$.  As in the symplectic case, the vanishing of $\Omega$
on ${\cal J}$ means that this form is symmetric and $J$-invariant.  It
is usual to suppose further that the normal bundle $TY/C$ has a
prescribed orientation, in which case it makes sense to require that
the Levi form be {\em positive} definite; in the negative case, we
speak of strict {\bf pseudoconcavity}.  Following standard terminology
in the symplectic case, we will call a strictly pseudoconvex complex
polarization a {\bf positive} polarization.

We note that the space of adapted complex structures on a symplectic
vector space, i.e. those for which the form $(x,y)\mapsto
\Omega(Jx,y)$ is positive definite and symmetric, is contractible.
Any two such almost complex structures are related by a transformation
which preserves $\Omega$ (which is therefore unitary); furthermore,
this transformation can be chosen in a ``natural'' way if one uses
the riemannian geometry of the symmetric space $Sp(2n-2)/U(n-1)$ to
select the  geodesic connecting the two structures and then lift it to
the symplectic group.

A CR structure is called {\bf embeddable} if there are enough CR
functions to realize $Y$ as the pseudoconvex boundary of a compact
normal (possibly singular) Stein domain $X_J$ (which is then uniquely
determined by $J$).  In dimension at least 5, all strictly
pseudoconvex CR structures are embeddable \cite{bo:integration} , but
in dimension 3 this is a real restriction.  The importance of
$X_J$ is that the smooth CR functions on $Y$ are precisely the
boundary values of holomorphic functions on $X_J$.  We refer to
\cite{ja:cr} for a general treatment of geometry and analysis on CR
manifolds.  

Epstein \cite{ep:relative} has shown that, if $J_1$ and $J_2$ are
embeddable CR structures on $Y$, then the orthogonal projection from
$H_{J_1}$ to $H_{J_2}$ is a Fredholm operator whose homotopy class is
independent of the choice of smooth measure on $Y$.\footnote{Actually,
the cited papers only prove this statement when $Y$ is 3-dimensional,
but the methods should extend to the general case.}  The
{\bf relative index} of $J_1$ and $J_2$ is thus finite in this situation.
Surprisingly, perhaps, the index is not always conserved under
deformations of $J_1$ and $J_2$.

For a positive polarization, the filling used for computing relative
indices will be taken to be the Stein domain mentioned above.   If 
the Stein domains $X_{J_1}$ and $X_{J_2}$ determined by a
pair of embeddable CR structures $J_1$ and $J_2$ on $Y$ are
nonsingular, these manifolds can be glued together along their common
boundary to form a closed manifold $X$.  Although the complex
structures on $X_{J_1}$ and $X_{J_2}$ do not match along $Y$, it is
possible, using the natural isomorphism between the vector bundle complex
structures mentioned above, to endow $X$ with a natural (up to
homotopy) stable almost complex structure and hence with a Dirac
operator $D^+$ which restricts away from a neighborhood of $Y$ to the
``rolled-up Dolbeault complexes'' (see \cite{gi:invariance}) on
$X_{J_1}$ and $X_{J_2}$.  Our gluing conjecture then states that the
relative index of $J_1$ and $J_2$ is equal to index of $D_+$.  We will
see in Section \ref{sec-singular}
we will see how to extend the conjecture to the singular case.

\section{Real polarizations}
\label{sec-real}

${\cal J}$ is a {\bf real polarization} if ${\cal J} = \overline {\cal
J} $.  This means that ${\cal J} $ is the complexification of the
tangent distribution of a foliation of $Y$ by legendrian submanifolds.
{\bf Fibrating} real polarizations are those for which this foliation
is a fibration.  Cosphere bundles foliated by their fibres are
examples of this type.  In fact, Pang \cite{pa:structure} proves that
these are the only examples with compact, simply connected leaves.
The quantum Hilbert space associated to $S^*M$ with its polarization
by fibres is just $L^{2}(M)$.  A filling in this special case is
constructed as follows.  Choose a riemannian (or finslerian) metric on
$M$, let $D^*M$ be the unit disc bundle in the cotangent bundle, and
identify the cosphere bundle with its boundary, so that the cotangent
disc bundle becomes the filling.

Given a contact transformation $\phi$ between
cosphere bundles $S^*M_1$ and $S^*M_2$, we may use it to {\em
identify} both bundles with a single contact manifold $Y$, which then
inherits a pair of real polarizations.  The quantum Hilbert spaces for
these polarizations are $L^2(M_1)$ and $L^2(M_2)$, but the operator
between then obtained by orthogonal projection in $L^2(Y)$ is not in
the class of Fourier integral operators $C(\phi)$ associated with
$\phi$ but is rather a {\bf Radon integral} operator associated with
the double fibration $M_1 \leftarrow Y \rightarrow M_2$.  This
operator, defined by pulling back by one fibration followed by
integration over the fibres of the other, is indeed a Fourier integral
operator, but its associated canonical relation is too big: it
contains at least the ``unoriented''
version of $\phi$ consisting of the graph of $\phi$ together with that
of $\xi\mapsto -\phi(-\xi)$, and is even larger except in ``clean''
cases. 

We should not, therefore, define the relative index of two real
polarizations to be the index of the orthogonal projection between
their quantum Hilbert spaces.  Instead, we must use an indirect
method, such as that described in the next section.  

\section{The Guillemin transform}
\label{sec-guillemin}
In order to realize Fourier integral operators as intertwining
operators between real polarizations, we follow an idea of Zelditch
and relate them through polarizations of CR type.  The groundwork for
this argument has been laid by Guillemin \cite{gu:toeplitz} in the
following way.  If $M$ is a compact manifold of dimension $n$, we
choose a real analytic structure on $M$ (which is essentially unique).
According to Grauert \cite{gr:levi}, $M$ can be
embedded as a totally real submanifold of a complex $n$-manifold
$M_{\complex}$ with strictly pseudoconvex boundary $Y$.  Like any
hypersurface in a complex manifold, $Y$ inherits a CR structure which
in this pseudoconvex case determines a contact structure on $Y$.  The
analysis of Guillemin shows that the {\bf Grauert tube}
$M_{\complex}$ can be identified
with a cotangent disc bundle $D^*M$ for some riemannian metric on $M$
in such a way that the contact structure $Y$ arising from
$M_{\complex}$ agrees with the one arising from the
identification of $Y$ with $S^*M$.  

$Y$ thus has two polarizations, one positive and one real.  We will
call these polarizations {\bf affiliated} with one another.  The
corresponding fillings are diffeomorphic, but one carries 
the structure of a Stein manifold while the other is 
symplectic.  Guillemin shows that the projection operator between the
quantum Hilbert spaces for these two polarizations (holomorphic
functions on $M_C$ in one case, and all functions on $M$ in the other)
is an elliptic Fourier integral operator with complex phase and hence
a Fredholm operator.  We will call this operator a {\bf Guillemin
transform} for $M$ and denote its index by $i_M$.  Guillemin shows
that this index is independent of all the choices made in its
construction and is therefore an invariant of the differentiable
manifold $M$.  Recently, Epstein and Melrose \cite{ep-me:shrinking}
have shown that this index is always zero.  In fact, they show that
the transform is an isomorphism for sufficiently small tubes.  (This
result can also be seen as verifying a special case of our gluing
conjecture, since the two fillings are topologically equivalent.)

Now for the idea of Zelditch \cite{ze:index}.  If $\phi$ is a contact
transformation between $S^*M_1$ and $S^*M_2$, we use it to identify
these two cosphere bundles with a common manifold $Y$ as before, but
now we consider {\em four} polarizations on $Y$.  In order, these are:
\begin{itemize}
\item
${\cal L}_1=$ the real polarization by fibres over $M_1$;
\item
${\cal J}_1=$ the positive polarization as the boundary of $M_{1,\complex}$;
\item
${\cal J}_2=$ the positive polarization as the boundary of $M_{2,\complex}$;
\item
${\cal L}_2=$ the real polarization by fibres over $M_2$.
\end{itemize}
Zelditch observes that the successive composition of the orthogonal
projections operators between the quantum Hilbert spaces of these
polarizations {\em is} a Fourier integral operator in the class
$C(\phi)$, so that its relative index can be computed as a relative index of
Epstein type between the two complex polarizations plus the difference
of Guillemin indices $i_{M_1}$ and $i_{M_2}$.  As we noted above,
the Guillemin indices are zero.
Thus, the index problem for Fourier integral operators is
reduced to the relative index problem for CR structures.

In general, to define the relative index between two polarizations, we
replace any which are real by affiliated positive polarizations.

\section{Extension to vector bundles}
\label{sec-extension}

The standard index theorems for pseudodifferential and Toeplitz
operators are most interesting when applied to operators on sections
of vector bundles rather than just on scalar functions.  The same
should be true for Fourier integral operators and their variants.  In
this section, we will propose a setup for an extension of our
conjectures to vector bundles, and we will see that the conjecture
reduces to known theorems in the pseudodifferential case.

Our starting data will now be a vector bundle $F$ over the contact
manifold $Y$, together with polarizations ${\cal J}_j$ of $Y$ corresponding
to fillings $X_j$.  In order to extend the vector bundle over the
fillings in an appropriate way, we need a condition of compatibility
with the polarizations.  In both the real and complex cases, the condition
will be ``constancy of the fibres along the leaves.''  

In the real case, where $X_j$ is a cotangent disc
bundle $D^* M_j$ and ${\cal J}_j$ is the polarization by fibres of the
cosphere bundle, the fibres of $F$ should be identical over all the
points of each fibre, which means that $Y$ should be the pullback to
$S^*M_j$ of a vector bundle $V_j$ over $M_j$.  In this case, we can
also pull back $V_j$ to the filling $X_j$ to give an extension of $F$
to a bundle whose fibres are constant along the leaves of the
polarization of the symplectic manifold $X_j$ by fibres of the
cotangent bundle.

In the complex case, we interpret ``constancy along the leaves of a
polarization'' as the existence of a flat connection along the
corresponding distribution.  When $J_j$ is a CR structure on $Y$, this
leads directly to the condition that the bundle $F$ should be a
{\bf holomorphic} vector bundle in the sense of Tanaka \cite{ta:differential}
(called an {\bf almost CR} vector bundle by Webster \cite{we:integrability}).
In this situation, we will further assume that $F$ extends to a
holomorphic vector bundle $E_j$ over the Stein filling $X_j$, and that
the CR sections of $F$ are the boundary values of holomorphic sections
of $E_j$.\footnote{It would be interesting to have verifiable
hypotheses guaranteeing the
existence of such extensions.}
The simplest example of this setup occurs when $F$ is a trivial
bundle, in which case we are simply dealing with $\complex ^N$-valued
functions which are CR on $Y$ and holomorphic on $X_j$.  

Once we have lifted the polarizations ${\cal J}_j$ on $Y$ to the
vector bundle $F$ as described above, we can identify a space of
smooth sections which are ``parallel in the direction of the
polarization,'' and then form their $L^2$ closure, using a volume
element on $Y$ and a hermitian structure on $F$, obtaining a space
which we will again call $H_j$.  The index of the orthogonal
projection from one space to the other is again well defined in many
cases and could be called the relative index of the two lifted
polarizations.  (When a polarization is real, we replace it by an
affiliated positive one before computing the index.)  As before, we
conjecture that this relative index is equal to the index of a Dirac
operator on the glued manifold $X$.  This time, the operator is a
twisted Dirac operator, obtained by tensoring with the vector bundle
over $X$ obtained by gluing the bundles $E_j$ by using their
identifications with $F$ over the common boundary $Y$.

We recover standard index theorems for Toeplitz and pseudodifferential 
operators by choosing the polarizations ${\cal J}_J$ (and
hence the fillings $X_j$) to be equal to one another, but by allowing
two different lifts of the polarizations to $F$.  For instance, if we
are given a bundle automorphism $\sigma$ of $F$, we can define one
lift to be the pullback of the second by $\sigma$.  In this case, if
$\pi$ denotes the orthogonal projection onto $H$ (which does
not depend on $j$ in this case), the operator which gives the relative
index $\pi \sigma\pi:H\mapsto H$.  When ${\cal J}$ is
positive, this operator is just the Toeplitz operator whose symbol is
$\sigma$, and our conjecture for the index reduces to the index
formula of Boutet de Monvel \cite{bo:index}.

When the polarizations are both real, with $X_j=D^*M_j$, $\sigma$ is
the symbol of a pseudodifferential operator $P$ between sections of
vector bundles $V_1$ and $V_2$ over $M$, The index of our glued
twisted Dirac operator is now the Atiyah-Singer topological index of
$P$, but the operator obtained from $\sigma$ by the projection process
described above is {\em not} $P$; rather, it is simply the
multiplication operator by the bundle map from $V_1$ to $V_2$ given by
integrating $\sigma$ over the fibres of the cosphere bundle $Y$.  To
get the operator $P$, we must use affiliated polarizations as
described in Section \ref{sec-guillemin} and use the results of
\cite{gu:toeplitz}.

\section{Singular fillings}
\label{sec-singular}

There are several ways to  approach the problem of
singular fillings.  One is to resolve the singularities and then add a
correction term to account for the nontrivial pseudoconvex (but no
longer Stein) filling.  We will present here an alternative approach
which appears to be more conceptual in nature.  It still uses
resolution of singularities, for the moment, but only to show that a
certain index is well defined, not to define it.

As usual, we consider polarized contact manifolds $Y$ of either of two
types--cosphere bundles and embeddable CR manifolds.  In the first
case, the filling will be the corresponding disk bundle in a
cotangent bundle; in the second, the filling will be the (possibly
singular) Stein domain having $Y$ as its strictly pseudoconvex
boundary.  

Let $X_1$ and $X_2$ be fillings of $Y$ corresponding to
polarizations 
${\cal J}_1=$ and ${\cal J}_2=$.  We may glue $X_1$ to $X_2$ along $Y$
to get a new object $X$, but the nature of $X$ depends on the nature
of $X_1$ and $X_2$.  If $X_1$ and $X_2$ are both either symplectic or
are nonsingular Stein varieties, they can be considered as almost
complex manifolds and hence $X$ becomes a stable almost complex
manifold.  The index of the corresponding Dirac operator is therefore
well defined.  

If either $X_1$ and $X_2$ is possibly singular, we will resort to the
following construction.  According to Theorem 8.1 of
\cite{le:algebraic} (see \cite{de-le-sh:algebraic} for related
results), 
each $X_j$ can be completed by adding a nonsingular complex manifold $Q_j$
with strictly {\em pseudoconcave} boundary $Y$ to make a (possibly
singular) projective variety $Z_j$.\footnote{I learned about this
result in a talk by G. Mati\'c on the paper \cite{li-ma:tight}, where I
also learned about gluing complex and symplectic manifolds along
contact boundaries!}  For such a variety, we {\em
define} the ``index of its Dirac operator'', denoted simply by $\indx
(Z_j)$ to be the Euler characteristic of its cohomology with values in
the sheaf ${\cal O}$ of germs of holomorphic functions.  This is a good
definition because, if $Z_j$ happens to be singular, this Euler
characteristic equals the Euler characteristic for the Dolbeault
cohomology on forms of type $(0,q)$, which is in turn equal to the
index of the Dirac operator given by the rolled-up Dolbeault complex.

If $Q_1$ and $Q_2$ were isomorphic, it would be reasonable to define
the relative index of $X_1$ and $X_2$ to be the difference of the
indices of the $Z_j$.  In general, account for the
difference  between $Q_1$ and $Q_2$ in the following
way.  Glue $Q_1$ and $Q_2$ along their common boundary $Y$ to form
a smooth manifold $Q$.  The complex structures on the pieces glue to
give a stable almost complex structure on $Y$ for which the natural
orientation agrees with that on $Q_2$ but is opposite to the
orientation of $Q_1$.  We now define the topological 
relative index of $X_1$ and
$X_2$ to be $\indx (Z_2) -\indx(Z_1) - \indx (Q)$, where the last
index is the index of the Dirac operator on $Q$ associated with its
almost complex structure.

Since the ``caps'' $Q_1$ and $Q_2$ are not unique, we have to check
that our relative index is well-defined.  This can be done by an
argument which we will not give here.  It uses the cobordism
invariance of the index and resolution of singularities.  (We hope
that the latter may be replaced by a localization argument for the
index of a singular variety.)

Our conjecture is that this expression $\indx (Z_2) -\indx(Z_1) -
\indx (Q)$ plays the role of the index of the object $X$ obtained by
gluing $X_1$ and $X_2$ along $Y$, and hence is equal to the relative
index of $X_1$ and $X_2$.  Using Riemann-Roch theory, it is not hard
to verify that the conjecture gives the correct relative index for the
pairs of $CR$ structures on a circle bundle over a Riemann surface of
genus 2 as considered in \cite{ep:relative}.

{\bf Remark} It would interesting to define the index of $X$ directly.
As a geometric object, $X$ can be thought of as consisting of two ends
which are (possibly singular) complex varieties, joined by a band on
which there is a stable almost complex structure.  The Dirac operator
of the band agrees on the overlap with the rolled up Dolbeault complex
on the smooth parts of the ends.

It is tempting to try to define the index of the glued object as the Euler
characteristic of an object in a derived category of sheaves on $X$,
obtained by gluing the sheaf $\cal O$ on the holomorphic ends to the
(very short) complex of sheaves given by the Dirac operator on the
band, using the techniques in \cite{ka-sc:sheaves}.  Unfortunately,
these sheaves are not quite quasi-isomorphic on the overlap of the two
regions---it is only the alternating sums of their cohomologies which
agree in some sense there.  Perhaps suitable holomorphic vector fields
near $Y$ could be used, in the spirit of \cite{at:vector}, to surmount
this problem.

\section{Holomorphic vs. Dirac indices: a proof strategy}
\label{sec-strategy}
Our strategy for proving the gluing conjecture for the relative index of CR
structures is to reduce the problem to related known results about
Dirac operators.  If $D^+$ is a Dirac operator between sections of
Clifford bundles $E^+$ and $E^-$ over a filling of the compact
manifold $Y$, then a famous result of Seeley \cite{se:singular} 
implies that the
orthogonal projection (the so-called {\bf Calderon projector}) 
from $L^2(Y)$ to the {\bf Cauchy data space} of boundary values of
solutions of $D^+ u = 0$ is a pseudodifferential operator of classical
type (i.e. with symbol an asymptotic sum of homogeneous terms) whose
principal symbol is a projection operator on the pullback of $E^+$ to
$S^*Y$.  Given a pair of such operators with Calderon projectors
having the same principal symbols, the orthogonal projection operator
between their Cauchy data spaces is shown to be a Fredholm operator by
Boo\ss-Bavnbek and Wojciechowski \cite{bo-wo:elliptic}, who prove the
following ``gluing theorem'' (originally conjectured by Bojarski) 
for the index of
this operator, which we call the {\bf relative index} of the two Dirac
operators.  (In general, it depends on the boundary isomorphism as
well as the operators.)

{\bf Theorem.}  {\em Let $D^+_1$ and $D^+_2$ be Dirac operators on compact
manifolds $X_1$ and $X_2$ having the common boundary $Y$, with
isomorphisms over $Y$ between the domain and range Clifford bundles,
such that their Calderon projectors have the same principal symbol
with respect to the domain isomorphism.  Then the relative index of
$D^+_1$ and $D^+_2$ is equal to the index of a Dirac operator on the
glued manifold $X=X_1 \cup _Y X_2$ obtained by gluing the bundles and
operators over $X_1$ and $X_2$ via the isomorphisms over $Y$.}

The Dirac operators to which we wish to apply the theorem above are
the Dolbeault-Dirac operators on the Stein fillings (assumed
nonsingular) $X_1$ and $X_2$ associated with a pair of positive polarizations
on the contact manifold $Y$.  More precisely, we assume that these
fillings are equipped with K\"ahler metrics (for instance those
obtained from embeddings in some $\complex^N$), and we consider on
each the operator $D^+ = \delbar+\delbar ^*:\Omega^{0,even}\rightarrow
\Omega^{0,odd}$ between the even and odd parts of the Dolbeault
resolution of the sheaf of holomorphic functions.  ``Rolling up'' the
Dolbeault complex by replacing its usual $\integers$ grading by a
$\integers _2$ grading has the result of replacing the rather delicate
Dirichlet problem for the $\delbar$ operator by a much more
robust problem, to which the gluing result above may be applied.  

The isomorphism over $Y$ between the domain and range bundles for 
$D^+_1$ and $D^+_2$ is obtained from an isomorphism between the
restrictions to $Y$ of the complex vector bundles $TX_1$ and $TX_2$.
This isomorphism is in turn obtained from the natural isomorphism between the 
two induced almost complex structures on the fixed contact
distribution $C$, as described in Section \ref{sec-complex} above.

The problem is now reduced to the following conjecture, in some sense
a relative version of the result in the compact case that the
dimension of the space of holomorphic sections of a line bundle
without higher cohomology is equal to the index of a rolled-up
(twisted) Dolbeault complex.  

{\bf Conjecture.}  {\em Let $X_1$ and $X_2$ be a nonsingular Stein fillings
of a contact manifold $Y$.  Then the relative index of $X_1$ and $X_2$
defined by the boundary values of their spaces of holomorphic
functions is equal to the relative index of the Dirac operators
$D^+_1$ and $D^+_2$.}

Some evidence in favor of this conjecture comes from the case where the
complex dimension of $X_j$ is 2.  In this case, the Cauchy data space
for the Dirac operator can be written as the direct sum (but not an
orthogonal one!) of the Cauchy data space for the holomorphic
functions and a subspace isomorphic to that for the harmonic forms of
type $(0,2)$.  The latter space is independent of the CR structure,
since the Dirichlet problem for the laplacian can be solved for any
Cauchy data.  Thus, in considering the relative indices for $X_1$ and
$X_2$, it ought to be possible to ``cancel'' the contributions coming
from the harmonic forms.

\end{document}